\documentclass[12pt]{article}
\usepackage{psfrag}
\usepackage{epsfig}
\usepackage{amsmath,amsthm,amsfonts,amssymb,amscd}
\theoremstyle{plain}
\newcommand{\diff}{\operatorname{Diff}}

\newcommand{\nat}{{\rm I\!N}}

\newcommand{\re}{{\flushleft{\bf Remark: }}}

\newcommand{\ra}{{\rightarrow}}

\def\cC{\mathcal{C}}
\def\cK   {\mathcal{K}}
\def\cF{\mathcal{F}}

\def\cU{\mathcal{U}}
\def\cD{\mathcal{D}}
\def\cP{\mathcal{P}}
\def\cV{\mathcal{V}}
\newtheorem{ttt}{Theorem}
\newtheorem{lem}{Lemma}[section]
\newtheorem{defi}[lem]{Definition}

\newtheorem{rem}[lem]{Remark}
\newtheorem{que}{Question}
\newtheorem{pro}[lem]{Proposition}
\newtheorem{coro}[lem]{Corollary}

\begin{document}
\title {Stably ergodic diffeomorphisms which are not partially hyperbolic }


\author{\sc Ali Tahzibi}
\date{February 01, 2002}
\maketitle
\begin{abstract}
 We show stable ergodicity of a class of conservative diffeomorphisms of $\mathbb{T}^n$
 which do not have any hyperbolic invariant subbundle. Moreover, the uniqueness of SRB
 (Sinai-Ruelle-Bowen) measure for non-conservative $C^1$ perturbations of such diffeomorphisms
 is verified. This class strictly contains non-partially hyperbolic robustly transitive
 diffeomorphisms constructed by Bonatti-Viana \cite{BoV00} and so we answer the question
 posed there on the stable ergodicity of such systems.
\end{abstract}
\section{Introduction}
One of the main aims of dynamical systems is to answer the following
questions:
\begin{enumerate}
\item Are the important topological or metric properties satisfied by majority of dynamical systems?
\item Under which conditions such properties persist after small perturbation of the system?
\end{enumerate}
%
%
%
%
%
%
Ergodicity is a basic feature of conservative dynamical systems
that yields the description of the average time spent by typical
orbits in different regions of the phase space. For
non-conservative systems the existence of SRB measures is a
natural candidate for the same purpose and they are defined as
follows.

 Let $M$ be a compact manifold and $f: M \ra M $. Given an
$f-$invariant Borel probability
 measure $\mu$, we call basin of $\mu$ the set $ B ( \mu )$ of $ x \in M $ such that:

$$  \lim_{ n \rightarrow \infty }  { 1 \over n } \sum_{ j = 0 }^{ n - 1 } \phi ( f^i ( x ) ) = \int \phi d \mu  \quad \text{for every } \quad  \phi \in C^0 ( M )   $$
and say that $ \mu$ is a physical or SRB (Sinai-Ruelle-Bowen) measure for $f$ if $ B ( \mu ) $ has
 positive Lebesgue measure .

A program proposed a few years ago by Palis~\cite{Pa00} contains a
conjecture  related to the first question:

Every system can be $C^r$ approximated, any $r \geq 1$, by one having finitely many SRB measures with their basins covering a full Lebesgue measure of the phase space.

By the above conjecture, we expect that for a``majority" of diffeomorphisms,
 the average time spent by typical orbits in different regions of the phase
 space is described by at most a finite number of measures.

In the same direction, in~\cite{ABV00} the authors show the
existence of finitely many SRB measures with basins covering a
full Lebesgue measure of the ambient manifold, for a large class
of partially hyperbolic systems and more generally for systems
displaying dominated splitting.

Let $M$ be a closed, compact riemannian manifold with volume form
$\omega$. A $C^2$-volume preserving diffeomorphism $ f : M \ra M $
is stably ergodic if there is a neighborhood $ \cU $ of $f$ in $
\diff_{\omega}^2 ( M ) $, the space of $ C^2$-volume preserving
diffeomorphisms of $M$, such that every $ g \in \cU $ is ergodic.

Considering the question (2), we want to know the necessary conditions to get stable ergodicity.
 First, Anosov in~\cite{An67} proved ergodicity of the Lebesgue measure for globally hyperbolic systems. Later, Pugh and Shub proved stable ergodicity of a large class of partially hyperbolic systems. The main condition to get ergodicity in these
results is ``accessibility": any two points of the phase space can be joined by a $C^1$-path consisting of consecutive segments, which are part of stable or unstable foliations
 (see \cite{BPSW} for a recent result in stable ergodicity).

  Recently F. Rodriguez~\cite{F01} showed the stable ergodicity of partially
   hyperbolic automorphisms
   of $\mathbb{T}^n$ for which the accessibility is not satisfied.

In this paper we show the stable ergodicity of an open set in
$\diff_{\omega}^1 ( \mathbb{T}^n )$ admitting no invariant
hyperbolic subbundle. In particular, we answer the question posed
in~\cite{BoV00} about stable ergodicity of the constructed
robustly transitive example there. The novelty of our work can be
explained as follows.

The existence of invariant foliations tangent to the hyperbolic subbundles of
 partially hyperbolic systems (\cite{HPS77}) is the main tool for proving the ergodi\-city
 of such systems.
 In the present work, no invariant hyperbolic subbundle is available.
  We have a dominated splitting and a non-uniform hyperbo\-licity property for this splitting
   which is explained in the Preliminary section.
  Moreover, we show the uniqueness of the SRB measures constructed in~\cite{ABV00},
  for non-conservative perturbations.

The class $\cV \subset \diff^1 ( \mathbb{T}^n )$ under
consideration consists of diffeomorphisms which
 are deformation of an Anosov diffeomorphism. To define $\cV$, let $f_0$ be a linear
  Anosov diffeomorphism of $n$-dimensional torus $\mathbb{T}^n$ (in fact, we need $f_0$ only
   to be an Anosov diffeomorphism on $ M = \mathbb{T}^n$ whose foliations lifted to $\mathbb{R}^n$
   are global graphs of $C^1$ functions). Denote by $ TM = E^s \oplus E^u $ the hyperbolic
    splitting for $f_0$ with dim $ ( E^s ) = s$, dim $( E^u ) = u$ and let $V = \bigcup V_i$
    be a finite union of small balls. We suppose that $f_0$ has at least a fixed point outside $V$.
     and say that $f \in \cV $, if it satisfies the following open conditions in $C^1$ topology :
%
%
%
%
%
%
%
\begin{enumerate}
\item $ TM $ admits a dominated decomposition and there exists small conti\-nuous cone fields
$C^{cu} , C^{cs}$ invariant for $Df $ and $D f^{-1}$ containing
respectively $ E^{u}$ and $E^{s}$
\item $ f $ is $C^1$ close to $f_0$ in the complement of $V$, i.e for $ x \notin V $ there is $\sigma < 1$:
$$ \| ( Df | T_x D^{cu}  )^{-1} \| < \sigma \quad \text{and} \quad \|   Df | T_x D^{cs} \| < \sigma$$
\item There exists some small $\delta_0 > 0$ such that for $x \in V$:
$$\| ( Df | T_x D^{cu}  )^{-1} \| <  1 + \delta_0 \quad \text{and} \quad    \|  ( Df | T_x D^{cs} \| < 1 + \delta_0  $$
where $D^{cu} , D^{cs}$ are disks tangent to $ C^{cu} $ and $ C^{cs}$
\end{enumerate}
%
%
%
%
%
\begin{ttt} \label{main}
Every $ f \in \cV \cap \diff_{\omega}^2 ( \mathbb{T}^n )$ is
stably ergodic.
\end{ttt}

For non-conservative diffeomorphisms in $\cV$ we prove the
uniqueness of SRB measures. An important property required in this
case, called ``volume hyperbolicity", is defined as follows.
\begin{defi} \label{flavio}
Let $f : M \ra M$ be a $C^1$ diffeomorphism and $TM = E^1 \oplus
E^2$; we say that this decomposition has volume hyperbolicity
property, if for some $ C > 0 $ and $\lambda < 1 $ :
$$ | det ( Df^n ( x ) | E^1 ) | \leq C \lambda^n ,  | det ( Df^{-n } (x ) | E^2) | \leq C \lambda^n  $$
\end{defi}
 \begin{ttt} \label{SRB}
Any $f \in \cV \cap \diff^2 ( \mathbb{T}^n )$ having volume
hyperbolicity property for $TM = E^{cs} \oplus E^{cu}$ has a
unique SRB measure with a full Lebesgue measure basin.
\end{ttt}
 In Theorem \ref{SRB}, the
volume hyperbolicity has the main role for proving non-uniform
hyperbolicity. Roughly speaking, by means of this property and a
good geometry of the invariant leaves of $f_0$, typical orbits do
not stay a long time in $V$.

In Section 3, we give an example of non-partially hyperbolic
diffeomorphisms which satisfy the hypothesis of Theorem~\ref{main}
and \ref{SRB}. Observe that although our class is not partially
hyperbolic, some weak form of hyperbolicity called ``dominated
splitting" exists. To justify this dominated splitting condition
we make the following comments:

\begin{enumerate}
\item
 Having a unique SRB measure with full support
 in a
 robust way requires some weak form of hyperbolicity.
 Namely if $ \cU $ is a $C^1$ open set of
 diffeomorphisms such that any $g \in \cU \cap\diff^2 (M)$ has an
 SRB measure $\mu$ with $supp (\mu) = M$
   then any $f \in U$ admits a dominated decomposition. (See Appendix A.)

\item
The persistence of positive measure sets of invariant tori due to Kolmogrov, Arnold, Moser, Herman
 and others shows that of course  some form of hyperbolicity is needed
 to get ergodicity. On the other hand, one hopes that stable ergodicity implies dominated splitting.
\end{enumerate}

In Section 2, we give some definitions which will be used in the
rest of the  paper and in Section 3 the example of robustly
transitive and non-partially hyperbolic diffeomorphisms of
$\mathbb{T}^n$ is constructed.

 In Sections 4 and 5, we analyse the geometry
of the basins of the SRB measures constructed in \cite{ABV00} for
systems with dominated splitting. It is shown that for each such
measure there exists some disk almost contained in the basin of it
and the radius of the disk is large enough to intersect the stable
manifold of a fixed point $q$ outside $V$. As the intersection of
$W^s ( q )$ with the mentioned disk in the basin of each SRB
measure is transversal, we can $C^1$ approximate these disks by
the $\lambda$-lemma.

After approximating the basins of two SRB measures the idea now, is to
apply some local accessibility argument.
In Sections 6 and 7, we prove the existence of local stable manifolds and
 absolute continuity of their holonomy for a positive measure subset of unstable manifold of $q$.  By this we prove that $ B ( \mu_i ) \cap B ( \mu_j ) \not= \emptyset$, for any two SRB measures. Then, the definition of the basin of SRB measures implies that $\mu_i = \mu_j$.

 It is worthwhile to emphasize that, in Pesin's theory for construction of local stable manifolds
  the set of regular points in the sense of Lyapunov plays a crucial role. By Oseledets' theorem
   they occupy a total probability subset of the ambient manifold, but in Theorem 1 we need to
    use these results for non-regular points. In Sections 6 and 7, we show that the coexistence
     of non-uniform hyperbolicity and a good control of the angle between the subbundles
     corresponding to non-uniform contraction and non-uniform expansion, enable us to construct
      stable manifolds and prove the absolute continuity of their holonomy.

Now the important point is that the union of the basin of SRB
measures constructed in~\cite{ABV00} contains a full Lebesgue
measure subset of the phase space. So, by our uniqueness result,
 for Lebesgue almost all $ x \in M $:
$$ \lim_{ n \rightarrow \infty }  { 1 \over n } \sum_{ j = 0 }^{ n - 1 } \phi ( f^i ( x ) ) = \int \phi d \mu  \quad \text{for every } \quad  \phi \in C^0 ( M ) .   $$
This is equivalent to the ergodicity of the Lebesgue measure for conservative diffeomorphisms.\\
{\bf Acknowledgements:} I would like to thank my advisor Professor Jacob Palis for his support and enormous encouragements during the preparation of this work. I also thank Professor Marcelo Viana, who suggested me the problem. The many conversations with him and Federico Rodriguez Hertz were very important to my work.
Finally, I wish to stress the remarkable research atmosphere at IMPA and to acknowledge financial supports of CNPq.

%
%
%
\section{ Preliminary}
We may consider some ways of relaxing uniform hyperbolicity, like:

\begin{itemize}
 \item non-uniform hyperbolicity
\item partial hyperbolicity
\item dominated splitting
\end{itemize}
{\bf Non-uniform hyperbolicity}:\\ This approach is due to
Pesin~\cite{Pe77} and it refers to
 diffeomorphisms with nonzero Lyapunov exponents in a full measure subset of phase space.
 Recall that $\lambda$ is a Lyapunov exponent at $x$ if
  $\lim_{n \ra \infty} \frac{1}{n} \log \| D_xf^n ( v ) \| = \lambda $ for some
   vector $v \in T_x M $. By Oseledets' theorem Lyapunov exponents exist for a total
    probability subset of $M$.\\\\
\
{\bf Dominated splitting}:\\ This approach is due to Ma\~{n}\'{e}
 and refers to
diffeomorphisms with a continuous decomposition of tangent bundle
of the phase space : $TM = E^{cs} \oplus E^{cu} $, with the
following property:
$$ \| Df | E_x^{cs} \| . \| Df^{-1} | E_{f ( x )}^{cu} \| \leq \lambda < 1 \quad \text{for all } \quad x \in M  $$
From here on just to emphasize the domination we write $TM =  E^{cs} \prec E^{cu}$.\\
Whenever we have a dominated splitting on $TM $, there are two
cone fields $C^{cu} , C^{cs} $ with the following properties:

  $$  C^{cu}_a ( x ) = \{ v_1 + v_2 \in E^{cs} \oplus E^{cu} ; \| v_1 \| \leq a \| v_2\| \} , \,\,\,\,\,\,\,\,\,\,\,\,\,  Df (C^{cu}_a ( x ) ) \subset C^{cu}_{\lambda a }  ( f ( x ) ) $$
  $$  C^{cs}_a ( x ) = \{ v_1 + v_2 \in E^{cs} \oplus E^{cu} ; \| v_2 \| \leq a \| v_1 \| \} , \, Df^{-1} (C^{cs}_a ( x ) ) \subset C^{cs}_{\lambda a } ( f^{-1} ( x ) )$$
A system for which $ TM = E^s \prec E^c \prec E^u $ is a dominated splitting and $E^s , E^u $ are respectively uniformly contracting and expanding (at least one of them is nontrivial) is called partially hyperbolic. If both uniform contracting and expanding subbundles exist, we call the diffeomorphism as ``{\it strongly partially hyperbolic}".\\\\
{\bf Key property: ``Non-uniformly hyperbolic" dominated splitting}\\
 To construct SRB measures for systems
  with a dominated splitting, by the methods in~\cite{ABV00} we need to verify ``non-uniform
   hyperbolicity" in a total Lebesgue measure set in the following sense. There is
   some $c_0 > 0 $ such that

\begin{itemize}
               \item
There exists a full Lebesgue measure set $H$ such that for $x \in H $ :
 \begin{equation} \label{jose}
  \limsup_{ n \rightarrow \infty } { { 1 \over n } \sum_{ j = 0 }^{ n - 1 } \log \| ( Df | E_{ f ^j ( x ) }^ { cu })^{ -1 } \| } \leq - c_0
\end{equation}
and also  :
\begin{equation} \label{alves}
\limsup_{ n \rightarrow \infty } { { 1 \over n } \sum_{ j = 0 }^{ n - 1 } \log \| Df | E_{ f ^j ( x ) }^ { cs } \| } \leq - c_0
 \end{equation}
\end{itemize}
We mention that the above conditions imply nonzero Lyapunov
exponents. Let us just mention that in the Pesin theory, some
invariant measure is fixed  and non-uniformly hyperbolic systems
refers to ones without zero Lyapunov exponent in a total measure
set. But, we are working with the Lebesgue measure which is not
invariant for non-conservative diffeomorphisms of $\cV$. In this
paper by non-uniform hyperbolicity we refer to the above
conditions.

To verify non-uniform hyperbolicity for the diffeomorophisms in
Theorems \ref{main} and \ref{SRB}, we use the volume hyperbolicity
property defined in the Introduction (Definition \ref{flavio}).

For non trivial examples of diffeomorphisms with  volume
hyperbolicity property we mention the following (see \cite{BDP}).
\begin{itemize}
\item {\bf Conservative systems}: Any $C^1$ conservative diffeomorphism with a
dominated splitting $TM = E^1 \prec E^2$ has the volume hyperbolicity property.
\end{itemize}
From this and the continuity of $det Df$, we conclude the
following corollary:
\begin{coro} \label{area}
For any $f \in \cV \cap \diff_{\omega}^2 ( \mathbb{T}^n )$, there
exists $\sigma_1 > 1$ and $C > 0$ such that $$| det ( Df^n ( x ) |
T_x ( D^{cs} ) ) | \leq C \sigma_1^{- n } , \quad | det ( Df^{-n }
(x ) | T_x ( D^{cu} ) ) | \leq C \sigma_1^{- n } ,$$ where $
D^{cs}, D^{cu}$ are disks tangent to $C^{cs}, C^{cu}.$
\end{coro}
%
%
%
%
%

 \section{ An example for Theorem ~\ref{main} }

 Here we give an example of systems that satisfy the hypothesis of Theorem~\ref{main}. The first non-partially hyperbolic and robustly transitive example is constructed in~\cite{BoV00} on
 $\mathbb{T}^4$. we apply their method and show that it works in larger dimensions.
 Let $ f_0$ be a linear Anosov diffeomorphism on the $ \mathbb{T}^n$ for which
$$ T_x ( \mathbb{T}^n ) = \mathbb{R}^n = E_1^s \prec E_2^s \cdots
\prec E_{n - 2}^s \prec E^u $$
where  dim $( E^u ) = 2$ and dim$ ( E_i^s ) = 1$.

 We may suppose that $ f_0$ has fixed points $ p_1 , p_2 , ... , p_{n - 2}$.
  Let $ V = \bigcup B ( p_i , \delta ) $ be a union of balls centered at
   $p_i$ and radius sufficiently  small $\delta > 0 .$
 The idea is to deform the Anosov diffeomorphism inside $V$, passing first through a
  flip bifurcation along $E_i^s \oplus E_{i + 1}^s$ inside $B_i =  B ( p_i , \delta )$
  and then other deformation (see fig 1), always composing with  discrete time  map of
   Hamiltonian vector fields to get volume preserving diffeomorphism. (For an example of such
   vector fields see ~\cite{BoV00})

More precisely, first we modify along stable direction $ E_i^s (
p_i ) \oplus E_{i + 1}^s ( p_i )$ for $1 \leq i \leq n - 3  $
until the index of $p_i$ drops one and two fixed points $ q_i ,
r_i$ are created. These new fixed points are of index $ n - 2 $.
In the next step composing with another Hamiltonian (two
dimensional), we mix the two contracting subbundles of $T_{q_i}M$
corresponding to $ E_i^s ( q_i )$ and $E_{i + 1}^s ( q_i )$. After
these deformations we have:
$$T_{q_i} M = E_1 \prec \cdots \prec E_{i - 1} \prec F_i \prec \cdots \prec E^u $$
where $F_i$ is two dimensional and corresponds to the complex eigenvalue.
Finally we do the same for $p_{ n - 2}$, but in the unstable direction of it.\\
\begin{figure}[t]
\begin{center}
\psfrag{p}{$p_i$}
\psfrag{Es1}{$ E^s_i ( p_i ) \oplus E^s_{i + 1} ( p_i)  $}
\psfrag{q1}{$q_i$}
\psfrag{q2}{$r_i$}
\psfrag{Fi}{$F_i$}
\epsfysize=1in
\centerline{\epsfbox{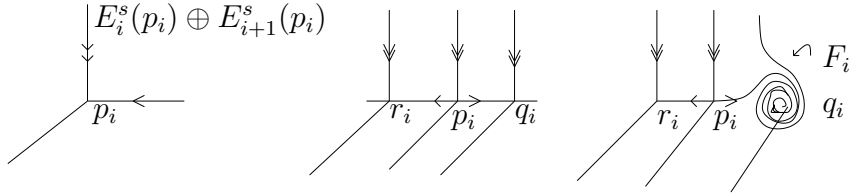}}
\caption{Deformation of Anosov}
\label{a}
\end{center}
\end{figure}

In this way we get an open set $ \tilde{\cV} $ in $ C^1$ topology
of diffeomorphisms satisfying the conditions 1-3 mentioned in the
introduction and:
\begin{itemize}
 \item There exist a hyperbolic fixed point $ q $ with stable index $ s = $ dimension
  of $E^s$ of the Anosov one (in the example is $ n - 2 )$, such that its stable manifold
   intersect any disk tangent to $ C^{cu}$ with radius more than $ \epsilon_0 $,
    for some small $\epsilon_0 > 0$. The similar thing for the unstable manifold
    and disks tangent to $C^{cs}$ happens. This is just because of the denseness
     of invariant leaves of the Anosov diffeomorphism $f_0$: Take a compact part
      of $ W^s ( q , f_0 ) $ to be $ \epsilon_0$ dense and taking $V$ small enough to garante
       permanence of this part during the deformations.
\end{itemize}
\begin{rem} \label{periodic}
Clearly the last item above is satisfied for $f \in \cV$ of
Theorems \ref{main} and \ref{SRB}, as $V$ is small enough.
\end{rem}
In what follows, we see that $f$ is robustly transitive but it is not partially hyperbolic.
 \begin{lem} \label{transitive}
$f \in \tilde{\cV}$ is robustly transitive.
\end{lem}
\begin{proof}  The proof goes as in $\mathbb{T}^4$ case in~\cite[Lemma 6.8]{BoV00} and we just
 remember the steps. The main idea to prove robust transitivity is to
  show the robust density of the stable and unstable manifold of an hyperbolic fixed point.
  We show the density of invariant manifolds of $q$ defined in Remark
  \ref{periodic} (see Proposition \ref{landa}).

  Let $U$ and $V$ be to open subsets. Using $\lambda$-Lemma and the
    density of invariant manifolds of $q$ we intersect some iterate of $U$ with $V$
   and get transitivity of $f$.
 \end{proof}
\begin{lem} $f \in \tilde{\cV}$ is not partially hyperbolic.\end{lem}
\begin{proof}
This is just because of the definition of partially hyperbolic systems: \\
$f$ is partially hyperbolic if $TM = E^s \oplus E^c \oplus E^u $
is a decomposition into continuous subbundles where at least two
of them are nonzero and $E^s$ and $E^u$ are respectively uniformly
contracting and expanding.
Suppose that $f$ is partially hyperbolic. First of all note that
by continuity of subbundles and the existence of a dense orbit by
lemma \ref{transitive}, the dimension of $E^s$ is constant.

 We claim that dim$ ( E^s ) = n - 2 $ and this gives a contradiction, because in $T_{p_i} M $
 there does not exist $n - 2 $ contracting invariant directions.
 To prove the claim observe that if  we suppose that dim $ ( E^s ) = j < n -2 $,
 then by the decomposition of $T_{q_j} M$ :
$$T_{q_j} M = E_1 \prec \cdots \prec E_{j - 1} \prec F_j \cdots \prec E^u .$$
 By definition, $E^s (q_j)$ must contain $E_1 \oplus \cdots \oplus E_{j - 1}$ and then as
 $F_j$ does not have any invariant subbundle we conclude that dim $ ( E^s ( q_j ) ) \geq j + 1$
 and this is a contradiction, because dim$ ( E^s ) = j$.\\
 By investigating $T_{p_{n -2}}M$, it is obvious that $f$ also can not have any continuous
  unstable subbundle.
\end{proof}
%
%
%
%
%
%
%
%
%
%
%
\section{$cu$-Gibbs measures}
 Gibbs measures in partially hyperbolic dynamical systems, as measures absolutely continuous
  along unstable foliation were constructed by Sinai and Pesin. (\cite{PS82})

For systems with only dominated splitting, in some cases we may call a probability
measure as $cu$-Gibbs, when its conditional measures with respect
to a measurable family of center-unstable disks  (tangent to
$C^{cu}$) is absolutely continuous with respect to the Lebesgue
measure of disks.

In \cite{ABV00}, SRB measures for the systems with dominated
decomposition having
 non-uniform hyperbolicity property and a technical called simultaneous hyperbolic times,
 are constructed.
  However in the Appendix(B) we show that, it is not necessary to verify such technical condition
   for constructing SRB measures of diffeomorphisms in $\cV$. The constructed SRB measures
    are in fact $cu$-Gibbs measures.
    Let us recall briefly the construction of $ cu$-Gibbs measures:
Fix a $ C^2$ disk  tangent to $ C^{cu}$ at every point of it and
intersecting $ H $ (the set of points having non-uniformly
hyperbolic behavior) in a positive Lebesgue measure where by
measure we refer the Lebesgue measure of the disk. Now consider
the sequence $\mu_n$ of averages of forward iterate of Lebesgue
measure restricted to such disk and then prove that a definite
fraction of each average corresponds to a measure $\nu_n$ which is
absolutely continuous with respect to Lebesgue measure along the
iterate of disk with uniformly bounded densities. Finally, show
that absolute continuity passes to $ \nu $, the limit of $\nu_n$.
More precisely:

\begin{pro} [\cite{ABV00}]
There exists a cylinder $\cC$ (a diffeomorphic image of product of two balls $ B^u , B^s $ of
 dimensions $ dim ( E^{cu} )$ and $dim ( E^{cs })$   in  $M$) and a family $ \cK_{\infty} $
 of disjoint disks  contained in $ \cC$ which are graph over $ B^u$ such that
\begin{enumerate}
\item The union of all the disks in $ \cK_{ \infty}$ has positive
$ \nu$ measure . \item The restriction of $\nu$ to that union has
absolutely continuous conditional measure along the disks in $
\cK_{\infty}$
 \end{enumerate}
\end{pro}
So we have $\mu = \nu + \eta$ where $\nu $ is absolutely continuous with a bounded away from zero Radon-Nikodym derivative along a $cu$-disks family. In this way we conclude that there exists disks $\gamma$ where $Leb_{\gamma}$-almost every point in $\gamma$ is regular and by absolute continuity of the stable manifolds ``for regular points", one gets a $\mu$ positive measure set in the same ergodic component. Normalizing the restriction of $\mu$ to the ergodic component above, we get an ergodic invariant probability measure $\mu^*$.

As the conditional measure of $\mu$ with respect to $ \cK_{\infty}$ is the sum of the conditional measures of $\nu$ and $\eta$ we conclude the following:
%
%
%

\begin{lem} \label{kinfinity}
There exists a disk $ D^{\infty}$ in $ \cK_{\infty}$, such that $Leb_{ D^{\infty}}$-almost every point of $D^{\infty}$ belong to the basin of $\mu^*$.
\end{lem}
 By Proposition 6.4 in \cite{ABV00}, $ M = \bigcup B ( \mu_i ) $ forgetting a negligible set,
  where $\mu_i $'s are $ cu$-Gibbs ergodic and SRB measures. By the above lemma
   we get the following corollary.
\begin{coro} \label{asli}
Let $f$ be as in Theorem~\ref{SRB}, then $ M = \bigcup B ( \mu_i ) $ (mod 0) where
 $\mu_i$ are ergodic SRB measures and for each $\mu_i$ there exists a disk $D_i^{\infty}$
  tangent to center-unstable cone field such that
  $D_i^{\infty} \subset B ( \mu_i)$ ($Leb_{(D_i^{\infty})}$-mod 0).
\end{coro}

In the next Section we prove that for $f \in \cV$, $B ( \mu_i  ) \cap B (\mu_j ) \neq \emptyset $ for all $i \not= j$, but as $  \mu_i$'s are ergodic so they are the same one.
\section{Uniqueness of $cu$-Gibbs measures}
{\bf Sketch of the proof of Theorems \ref{main} and \ref{SRB}}

\begin{figure}[t]
\psfrag{q}{$q$}
 \psfrag{Wu}{$W^u ( q )$} \psfrag{Ws}{$W^s ( q )$}
\psfrag{D1}{$D_i^{\infty}$} \psfrag{D2}{$D_j^{\infty}$}
\epsfysize=2in \centerline{\epsfbox{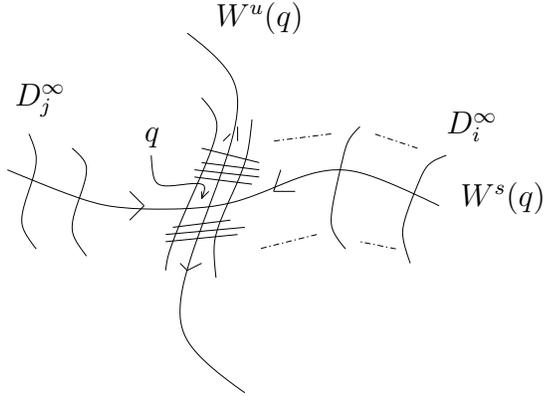}}
\caption{Intersecting basins via local stable manifolds}
\label{lambda}
\end{figure}

In order to show the uniqueness of the $cu$-Gibbs measures, we
prove that their basins have non empty intersection. For this, we
use Remark \ref{periodic} for the diffeomorphisms in $\cV$ and
Proposition \ref{landa} below to approximate the basins of SRB
measures. Then by means of local stable manifolds intersect the
basins corresponding to the different measures. Observe that any
two points in some local stable manifold belong to the basin of
the same measure. Let $q$ be as in Remark~\ref{periodic}.
\begin{pro} \label{landa}
The stable manifold of $ q , W^s ( q ( f ) , f )$ is dense and intersects
transversally each $ D_i^{\infty}$.
\end{pro}
\begin{rem}
This intersection is a crucial part of the proof of ergodicity.
Let's mention that just denseness does not imply intersection with
$D_i^{\infty}$.
\end{rem}
\begin{proof}
To prove the Proposition \ref{landa} we claim that some iterate of
$ D_i^{\infty} $ contains a disk tangent to center-unstable
conefield with radius more than $\epsilon_0$ which also almost
every point in it belong to $B ( \mu_i )$. This proves the
Proposition because of $\epsilon_0$ denseness of the $W^s ( q ) $,
see Remark \ref{periodic}. we prove the claim as following:

Consider a lift $\tilde{f} : \mathbb{R}^n \ra \, \mathbb{R}^n$ of
$f$, and $\pi_u$ as the projection along stable
 foliation of $f_0 $ (the Anosov one) from $ \mathbb{R}^n $ to
$\mathbb{R}^u$.
As $D_i^{\infty}$ is tangent to conefield $C^{cu}$ at each point
of it, we consider a global graph $\Gamma $ (the graph of a $C^1$
function $\gamma : \mathbb{R}^u \ra \, \mathbb{R}^s  , \| D \gamma
\| \leq \epsilon $ (angle of conefield), which contains
$D_i^{\infty}$.
 Now consider the iterates $\Gamma_n := \tilde{f}^n ( \Gamma ) $ which all of them
  are graph of $C^1$ functions with small derivative, this is because $f^n ( \Gamma ) $
   is a proper embedding of $\mathbb{R}^u$ in $\mathbb{R}^n$ whose tangent
    space at every point is
   in $C^{cu}$ and $C^{cu}$ is forward invariant.

 Now as $ Df $ expands area of disks on center unstable direction by arguments
 of ~\cite[Lemma 6.8]{BoV00} there exists some point $ x_0$ in $ f^{n_0} ( D_i^{\infty} )$
 such that its positive orbit never intersects $V$, so any small disk in $ \Gamma_{n_0} $
 around $ x_0 $ will have some iterate containing a disk with radius at least $ \epsilon_0 $.
  (See Remark \ref{periodic} for $\epsilon_0$). In this way we can
  show the denseness of $W^s (q)$. If $U $ is any
  open set just consider a  center-unstable disk $D$, in the intersection of
  $U$ and an unstable leaf of $f_0$ and argue as above
  substituting $D_i^{\infty}$ by $D$. The density of $Wû (q)$
  comes out by the similar method.

Now observe that because of the invariance of continuous cone field $C^{cs}$, the global
 stable manifold of $q$ is tangent to $C^{cs}$ at any point and consequently the intersection of
  $W^s (q)$ and $D_i^{\infty} $ is transversal.
\end{proof}
Using $\lambda$-lemma, for $ n$ large enough $f^n (D_i^{\infty} )$
and $ W^u ( q )$ are
 $ C^1 $ near enough. On the other hand, in Section \ref{non-uniform} we show that almost every
  point of $ W^u ( q ) $ have a local stable manifold.
This implies that there exists $S \subset W^u ( q )$ with $Leb ( S ) > 0 $ such that
 for all $x \in S$ the size
 of $W_{loc}^s (x) $ is uniformly bounded away from zero and $W^s_{loc}
 (x)$ intersects $f^n ( \cD_i^{\infty} )$ for $n$ large enough.
We need an absolute continuity property proved in Section
\ref{absolute continuity} to conclude the following:

 $$Leb_{ f^n ( \cD_i^{\infty} ) } \bigl( \bigcup_{ x \in S} W_{loc}^{s} ( x ) \cap B (  \mu_i ) \cap f^n ( \cD_i^{\infty} ) \bigr) > 0  $$
We would get the same thing for
 $ \mu_j$ and this enables us to find at least two points $ x , y $ respectively in  $ B ( \mu_i )$ and $ B ( \mu_j ) $  such that they
 are in the local stable manifold of the same point in $ S $ (see fig
 ~\ref{lambda}). This means

  $$ \lim_{ n \rightarrow \infty }  { 1 \over n } \sum_{ j = 0 }^{ n - 1 } \phi ( f^i ( x ) = \lim_{ n \rightarrow \infty } { 1 \over n } \sum_{ j = 0 }^{ n - 1 } \phi ( f^i ( y )  \quad \text{for every } \quad  \phi \in C^0 ( M ),  $$
and consequently $ B ( \mu_i ) \cap B ( \mu_j ) \ne \emptyset $
which implies $\mu_i = \mu_j $.
 We have proved that the decomposition of $ \mathbb{T}^n $ (mod 0) by the basin of SRB measures contains a unique element (mod 0) or there exists just one SRB measure whose basin has full Lebesgue measure.\\
%
%
%
If $f$ preserves the Lebesgue measure, by dominated splitting the
volume hyperbolicity is satisfied (see Preliminary). So by
Theorem~\ref{SRB} for almost all points
             $$ \lim_{ n \rightarrow \infty }  { 1 \over n } \sum_{ j = 0 }^{ n - 1 }
                     \phi ( f^i ( x ) = \int \phi d \mu
            \quad \text{for every } \quad  \phi \in C^0 ( M ) $$
 and immediately we have ergodicity of Lebesgue measure, completing the
 proof of the Theorem~\ref{main}.
\re In Theorem~\ref{SRB} we prove the uniqueness of the SRB
measures. The unique SRB measure is absolutely continuous along
disks which are unstable manifolds corresponding
  to positive Lyapunov exponents. By~\cite{LY85} one has the
  following:
$$ h_{\mu} ( f ) = \sum \lambda_i^{+} \quad \text{where}
\quad \lambda_i^{+} = max \{ 0 , \lambda_i \}, $$ where
$\lambda_i$ are the Lyapunov exponents of the ergodic measure
$\mu$. In fact, as the basin of the physical measure constructed
in Theorem~\ref{SRB} occupies
 a total Lebesgue measure set of manifold, it will be the unique measure among the ergodic measures
  with nonzero Lyapunov exponents which satisfy the Pesin's formula.

   we observe that with
   the same method with which we have proved the uniqueness of SRB measures, one also can show
    that $\mu $ is the unique ergodic measure satisfying Pesin's equality and having
     $u ( = $dim $ E^u$) positive Lyapunov exponents. Then by ergodic decomposition
     theorem it is the unique invariant probability with the mentioned
     properties. So, the following question is interesting:
  \begin{que}
  Does any $f$ as in Theorem~\ref{SRB}
  have only one measure satisfying the Pesin's equality?
  \end{que}
%
%
%
%
%
%
%
\section{ Non-uniform hyperbolicity } \label{non-uniform}

In this chapter we show a non-uniformly hyperbolic behavior  for a
full Lebesgue measure subset on
 the unstable manifold of the persistent hyperbolic fixed point $q$.
  Then we construct local stable manifold for the point of this subset.

To prove a non-uniform hyperbolic behavior, we will ``follow the
orbit of points" and observe that they spend a definitive part of
their time, out of the perturbation region and conclude that they
``remember hyperbolicity of the Anosov one". More precisely let
$W$ be a $u-$dimensional submanifold of $\mathbb{T}^n$ and $\pi$
the natural projection from $\mathbb{R}^n$ to $\mathbb{T}^n$. We
call $W$ dynamically flat according to the following definition.
\begin{defi}
$W$ is dynamically flat if for $\widetilde {W_n}$, any lift of
$f^n ( W ) $ to $ \mathbb{R}^n$, $ Leb (\widetilde {W_n} \cap K )
\leq C $ where $K$ is any unit cube in $\mathbb{R}^n$ and $C$ is a
constant depending only on $f$.
\end{defi}
\begin{lem}
             $W^u ( q )$ is dynamically flat.
\end{lem}
\begin{proof}
 Consider $\cF_0 ( q )$ the leaf of unstable foliation of $f_0$
 which passes through $ q $ and let
$ \cF_n = f^n ( \cF_0 ( q ) ) $. As $ \cF_0$ is a leaf of a linear
Anosov diffeomorphism, any lift of it to $ \mathbb{R}^n$
 will be a  $u$-affine subspace and is a proper image of
 $ \mathbb{R}^u$ to $ \mathbb{R}^n$. By invariance of the thin conefield
 $ C^{cu}$, we conclude that the tangent space of any lift of $
 \cF_n$, which we call also $\cF_n$,
 at every point is in $C^{cu}$ and it is also proper image of
 $ \mathbb{R}^u$. In this way for any unit cube $K$, $\cF_n \cap K$ can be seen as the graph
 of a $ C^1$ function with $u$-dimensional base of the cube as its domain.
  This function has an small
 norm of derivative which is independent
 of cube $K$ and $n$, this is because its graph is tangent
  to $C^{cu}$. So $\cF_n \cap K$ has a uniformly
 bounded area (with respect to Lebesgue measure of $\cF_n$)
  and this is what we want,
because the intersection of the unstable manifold of $ q $
 with $K$ is contained in the limit of $ \cF_n \cap K $. \\
\end{proof}
\begin{pro} \label{almost}
 Let $W$ be a dynamically flat submanifold and $f$ satisfying the hypothesis
  of Theorem~\ref{SRB}, then every small disk in $W$ contains a  Lebesgue total
   measure (Lebesgue measure of $W$) subset for which:
$$ \limsup_{ n \rightarrow \infty } { { 1 \over n } \sum_{ j = 0 }^{ n - 1 } \log \|
 ( Df | E_{ f ^j ( x ) }^ { cs }) \| } \leq - c_0   \hspace{.3 cm} c_0 > 0 $$
\end{pro}
\begin{proof}
Here we use the same argument of ~\cite{ABV00} and show that:

\begin{lem} ~\label{viana}
 There exists $ \epsilon > 0 $ and a total Lebesgue measure subset of any small disk $D$ in $W$, such that $ \# \{ 0 \leq j < n : f^j ( x ) \notin V \} \geq \epsilon n $ for every large n . \end{lem}
\begin{proof}
 we choose a partition in domains $ B_1 , B_2 \cdots B_{p + 1} = V $
 of $ \mathbb{T}^n$
such that there exists $ K_i , L_i  $ with $ B_i \in \pi ( K_i ) $
and $ f ( B_i ) \in \pi ( L_i ) $ where $K_i$ , $L_i$ are a finite
open cubes in $\mathbb{R}^n$)
 and estimate the Lebesgue measure of the sets $ [ \underline i ]$'s
where $ \underline i $ is an array $ ( i_0 , i_1 , ... , i_{n - 1} ) $
 and $ [ \underline i ]$ is defined as points in  $ D $ such that
 $ f^j ( x ) \in   B_{ i_j}$ for $ 0 \leq j < n $. In fact, we prove
  the following lemma. Let $\sigma_1$ be as in Corollary~\ref{area} then:
\begin{lem}
 $Leb ( [ \underline i ] ) \leq C \sigma_1^ { - n }$ (where $ C $ is a constant
  depending only to $f$)
\end{lem}

\begin{proof} By the choice of $B_i$ and induction we have that
 $ f^j ( [ \underline i ] ) \in \pi ( \widetilde W_n \cap L_{i_{ j - 1 }} ) $,
  where $\widetilde {W_n} $ is a lift of $f^n ( W )$ to $ \mathbb{R}^n$.\\
To conclude lemma we use area expanding (Corollary\ref{area}
 property along disks tangent to center unstable conefield and the fact that
  intersection of $\widetilde {W_n}$with a unit cube has a uniformly bounded volume.
By induction
   $$ Leb ( [ \underline i ] ) \leq \sigma_1^{ - n } Leb ( f^n ( [ \underline i ] ) \leq \sigma_1^{ - n } Leb ( \widetilde W  \cap L_{i_{ n - 1 }} ) \leq C \sigma_1^{ - n }$$
\end{proof}
Now we show how to conclude Lemma \ref{viana}. Let $ g (
\underline i ) $ be the number of values $ 0 \leq j \leq n - 1 $
for which $ i_j \leq p$. We note that the total number of arrays
with $ g ( \underline i ) \leq \epsilon n  $ is bounded by
$$ \sum_{k \leq \epsilon n } { n \choose k } p^k \leq \sum_{k \leq \epsilon n }
 { n \choose k } p^{ \epsilon n }$$
 and applying Stirling's formula gives that it is bounded by
  $e^{\beta_0 n} p^{\epsilon n}$ ($ \beta_0 $ goes to zero as
  $ \epsilon$ goes to zero).
So, the union of the sets $ [ \underline i ]$
 for which $ g ( \underline i ) \leq \epsilon n $ has Lebesgue
  measure less than
  $ C e^{\beta_0 n } p^{\epsilon n } \sigma_1^{ - n }$.
  Choosing $ \epsilon$ small enough such that
      $ e^{\beta_0 } p^{\epsilon} < \sigma_1$, we are in
       the setting of Borel-Cantelli lemma and conclude Lemma ~\ref{viana}.
\end{proof}
 By this lemma the Proposition \ref{almost} is proved just taking $ c_0 = - \log ( \sigma^{\epsilon} ( 1 + \delta_0 )^{ 1 - \epsilon} )$ and $\delta_0$ small enough.
\end{proof}
\begin{coro} \label{unstablemanifold}
 Almost all points of local unstable manifold of $q$ satisfy non-uniform hyperbolicity
  property.
\end{coro}
For any $x$ satisfying the conclusion of Proposition \ref{almost},
there exists $N ( x )$ such that for $ n \geq N $
$$\prod_{i = 0}^{n - 1} \| Df _{| E^{cs} ( f^i ( x ) ) } \| \leq
\lambda^n $$
  we remember that $ \lambda = \sigma^{\epsilon} ( 1 + \delta_0 )^{ 1 - \epsilon}$
   which is less than one if $\delta_0$ is small enough.
\begin{coro} \label{s}
There exists a positive Lebesgue measure subset $ S \subset W^u ( q ) $, $ N \in \rm I \!N$ and $ \lambda < 1$ such that
$ \forall x \in S $:
$$ \forall n > N \quad  \quad \prod_{i = 0}^{n - 1} \| Df _{| E^{cs} (  f^i ( x ) ) } \| \leq \lambda^n $$
\end{coro}
The points of $S$ are not necessarily regular in the sense of
Lyapunov. We can not use Pesin theory directly for the existence
of invariant manifolds and absolute continuity of their holonomy.
By dominated splitting and non-uniform hyperbolicity  as above we
can construct local stable manifolds.

 \begin{pro} Every point of $S$ has an stable manifold, whose size is bounded
  away from zero.
 \end{pro}
 \begin{proof}
We can construct local invariant disks using only domination property but in
 general case we do not know whether these disks are stable manifolds or not.
  For $f \in \cV$ by Corollary~\ref{s}, we are able to prove that the disks passing through
    the point of $ S $ are stable manifolds, i.e  $ d ( f^n ( x ) , f^n ( y ) ) \ra \,
 0 $ exponentially fast, for $ y \in W_{ loc }^{cs} ( x )$.\\
Denote $ Emb ( D^u , M ) $ the space of $ C^1$ embeddings from $ D^u$ to $ M $
 endowed with the $ C^1$ topology where $ D^u$ is the  $u$-dimensional ball of radius one.

Using notation of ~\cite{HPS77}, $ M $  is an immediate relative
pseudo hyperbolic set for $ f $ if there exists a continuous
function $ \rho $ such that:

 \begin{equation}\label{pugh}
  \| Df _{| E^{cs} ( x )} \| < \rho ( x ) < m (  Df _{| E^{cu} ( x )} )   \end{equation}
 where $ m ( T ) =  \| T^{-1} \|^{-1}$

%
%
%
 In our case, dominated splitting and compactness of $ \mathbb{T}^n $ imply
  relative pseudo hyperbolicity. We deduce that, there exist continuou sections
  \begin{itemize}
  \item $ \phi^u : M \ra Emb ( D^u , M ) $
  \item $ \phi^s : M \ra Emb ( D^s , M ) $
  \end{itemize}
such that $  W_{ \epsilon }^{cs} ( x ) := \phi^s ( x ) D_{ \epsilon }^s ,  W_{ \epsilon }^{cu} ( x ) := \phi^s ( u ) D_{ \epsilon }^u $, have the following properties:
\begin{enumerate}
 \item
 \begin{itemize}
 \item $ T_x W_{ 1 }^{cs} ( x ) = E^{cs} ( x ) $
 \item $ T_x W_{ 1 }^{cu} ( x ) = E^{cu} ( x ) $
 \end{itemize}
 \item Local invariance property . for all $ 0 < \epsilon_1 < 1 $ there is $ 0 < \epsilon_2 < 1$
 \begin{itemize}
 \item $ f ( W_{ \epsilon_2  }^{cs} ( x ) ) \subset W_{ \epsilon_1  }^{cs} ( f ( x ) ) $
 \item $ f^{-1} ( W_{ \epsilon_2  }^{cu} ( x ) ) \subset W_{ \epsilon_1  }^{cu} ( f^{-1} ( x ) ) $
 \end{itemize}
 \end{enumerate}

 Given any $ c $ we can take $ \epsilon_1 $ such that :
 \begin{equation}\label{continuity}
 1 - c < \frac  {  \| Df _{| T_y W^{cs} ( x ) } \| }{  \| Df _{| E^{cs} ( x )} \| } < 1 + c
  \hspace{.3 cm} \quad \text{when} \quad  \hspace{.3 cm}  d ( x , y ) < \epsilon_1 \,
       , y \in W_{ \epsilon_1 }^{cs} ( x )
\end{equation}
We can take this $ \epsilon_1$ uniformly in $x$ as $ M $ is compact and the section is continuous with image in embeddings endowed with $C^1$ topology.
 Choosing $ \epsilon_2$ such that $ f^i ( W^{cs}_{\epsilon_2} ( x ) ) \subset W^{cs}_{\epsilon_1} ( f^i ( x ) ) $ for all $ 0 \leq i \leq N $
    we show that
                   $ \forall n \in \nat ,\,\, d ( f^n ( x ) , f^n ( y ) ) \leq \epsilon_1 $.
 In fact, $d ( f^n ( x ) , f^n ( y ) )$ goes to zero as $n$ goes to infinity. We prove it by induction; let us define:

\begin{equation} \label{landabar}
 \overline{\lambda} : = ( 1 + c ) \lambda
\end{equation}
 and $ c $ is adjusted such that $ \overline{\lambda} < 1 $. As
 $d ( f^{i} ( x ) , f^{i} ( y ) ) \leq \epsilon_1 $ for $0 \leq i \leq N + k -1$
  we have:

$$ d ( f^{N + k} ( x ) , f^{ N + k} ( y ) ) \leq ( 1 + c ) \| Df_{| E^{cs} ( f^{ N + k - 1} ( x ) )} \| d ( f^{N + k - 1} ( x ) , f^{ N + k - 1} ( y ) ) \leq $$
$$ \leq \prod_{i = 0}^{ N + k - 2}  \| Df _{| T_{ z_i } W_{ \epsilon_1 }^{s} ( f^{i} ( x ) ) } \| \, \| Df_{| E^{cs} ( f^{ N + k - 1} ( x ) )} \| \, d ( x , y )$$
$$ \leq ( 1 + c )^{ n +k } \prod_{i = 0}^{ N + k - 1}  \| Df _{| E^{cs} ( f^i ( x ) )} \| d ( x , y ) \leq \overline{\lambda}^{ n + k } d ( x , y )$$
where  $z_i \in  W_{ \epsilon_1 }^{cs} ( f^{i} ( x ) ) $; this is all by Mean Value Theorem.
\end{proof}
\section{Absolute continuity} \label{absolute continuity}
\footnote{I thank Krerley Irraciel for useful discussions on this
section. } In this Section we prove that the holonomy map by the
local stable manifolds constructed  on $S$  is absolutely
continuous:
\begin{ttt} \label{absolute}
For large $n$, holonomy map from $ S \subset W^u_{local} ( q ( f )
, f ) $ to
  $ f^ n ( \cD_{\infty}^i ) $ is absolutely continuou i.e it sends
   the nonzero measure subset of $ S $ to a nonzero measure subset
    of $ f^ n ( \cD_{\infty}^i ) $.\end{ttt}

Let us mention that holonomy map $ h $ is defined on whole $S$ for
large $n$. From now on we call
 its inverse by $\pi$ which is holonomy along stable manifolds from
 $f^ n ( \cD_{\infty}^i ) $ to $ W^u ( q )$. We are going to prove that
 if $ B $ is a measure zero set in $ h ( S ) \subset f^ n ( \cD_{\infty}^i ) $
  then $ Leb ( \pi ( B ) ) = 0 $ and then conclude that $  Leb ( h ( S ) ) \neq 0 $.
   For this, it is enough to show that for every disk $ D \subset f^ n ( \cD_{\infty}^i ) $
 with center in $h ( S )$, the holonomy $ \pi$ from $ D$ to $ W^{u}( q )$
  does not increase measures, more than a constant which is uniform for all
   such disks :
      $$Leb ( \pi ( D ) ) < K  Leb ( D )$$
 because for any measurable set $ B $ with zero measure, we can cover
  it by a family of disks $ \cD $  such that $ \sum_{ D \in \cD } m ( D ) $
  is arbitrary close to zero. As $ Leb ( \pi ( D ) ) \leq  K Leb ( D ) $,
   we conclude that $ Leb ( \pi ( B ) ) = 0 $. From now on $S'$ represents $h ( S )$.\\
The proof of this absolute continuity result goes in the same
spirit of ~\cite{PS89}.
 The difference is that here the points for which we construct stable manifolds
 are not necessarily regular.
  We see that a nonuniform hyperbolicity and a good control
   on the angles of two invariant subbundles is enough to get an absolute
    continuity result. A short sketch of the proof is as follows.

To compare $Leb ( D )$ and $Leb ( \pi ( D ) )$, we iterate
sufficiently such that $ f^n ( D ) $ and $ f^n ( \pi ( D ) ) $
``become near enough". But after such iteration, $f^n ( D )$ may
have an strange shape, so in \ref{covering} we consider a covering
of $ f^n ( S' ) \cap f^n ( D ) $ by
  $  B_i := B ( a_n , f^n ( x_i ) ) $ (ball of radius $ a_n $ with center
   $ f^n ( x_i ) $) where $ x_i $ is in $ S' $ and $ a_n $ is much larger
   than $d ( x_i , \pi_n ( x_i ) )$ where $ \pi_n $ is defined naturally
    by $ \pi_n = f^n \circ \pi \circ f^{-n} .$

By the specific choice of $a_n$, in \ref{comparing} we show that
$Leb ( B_i ) \approx Leb ( \pi_n ( B_i ) )$.
 Indeed, the dominated splitting of the tangent bundle allows
 us to choose them in such a good way.
Finally, in \ref{distortion} we prove some distortion results and come back to
 compare the volume of $D$ and $\pi ( D )$.\\
%
%
%
%
%
%
%
%
\subsection{ Some general statements}
Let us fix some notations and definitions:
\begin{itemize}
\item $ d_1 ( $resp.$ d_2 ) := $ restriction of the riemannian
metric of manifold to
 $ f^n ( D ) $ $ ($resp. $ W^u ( q ) ) $
\item $ d_s := $ Intrinsic metric of stable manifolds
\item $ d $ := Riemannian metric of the manifold $M$
\item $ a \preceq b $ means $ a \leq kb $ for a uniform $ k > 0 $  $a , b \in 
\mathbb{R}  $
\item $a \approx b$ means that $k^{-1} < \frac{a}{b} < k$ for a uniform $k > 0$
\end{itemize}
\begin{defi} \label{angel}
If $ E , F $ are two subspaces of the same dimension in $
\mathbb{R}^{n}$, we define the angle between them $ \measuredangle
( E , F ) $, as the norm of the following linear operator : $$ L :
E \ra E^{\perp} \quad \text{such that} \quad  \, Graph ( L ):= \{
( v , L ( v ) ) , v \in E \} = F $$
\end{defi}

\begin{defi}A thin cone $ C_{\epsilon} ( E )$with angle $\epsilon$ around $E$ is defined
as subspaces $ S $ s.t $\measuredangle ( S , E ) \leq \epsilon .$
\end{defi}
By the definition of cones, it is easy to see that:
\begin{lem} \label{dominada}
If $ C^{cu}_{\epsilon} , C^{cs}_{\epsilon}$ are two conefield
which contain $ E^{cu}, E^{cs}$ (the subbundles of dominated
splitting) then $ Df_x C^{cu}_{\epsilon} ( x ) \subset C^{
cu}_{\lambda \epsilon} ( f (x ) ) $,  for some $ 0 < \lambda < 1 $
or in other words the angle will decrease exponentially.
 \end{lem}
\begin{proof}
 Take $S \in C^{cu} ( x ) $ and $v \in S$. By definition $v = v_1 \oplus v_2 $
where $ v_1 \in E^{cu} , v_2 \in E^{cs}$
and by dominated splitting (see Preliminary Section):
          $$ \frac{\| Df_x ( v_2 ) \|}{ \| Df_x ( v_1 ) \|} \leq
             \lambda \frac{\| v_2 \|}{\| v_1 \|}  $$
and this means that $ \measuredangle ( Df ( S ) , E^{cu} ( f ( x ) ) ) \leq \lambda  \, \measuredangle ( S , E^{cu} ( x ) ) $ by definition \ref{angel}.
\end{proof}
Let us state a lemma that gives us some good relations between $
d_1 , d_2 $ and $d.$
\begin{lem} ~\label{fasele}
If  $\mathbb{R}^n = S \oplus U ( U = S^{\perp} ) $ and $ h $ is a
$ C^1$ function from $ B ( 0 , \delta ) \subset S$ (ball of radius
$\delta$) to $ F $ where $ F $ is in a small cone $ C_{\epsilon} (
U ) $.  Suppose that $ T_x ( graph ( h ) ) \subset C_{\epsilon} (
S )  , \forall x \in graph ( h ) $ then:
\begin{itemize}
\item $ d_h ( z , 0 ) \leq C ( \epsilon ) d ( z , 0 ) ,\hspace{0.2 in} d_h : $distance on graph$ ( h )  $
\item $ Leb ( graph ( h ) ) \leq C (\epsilon) Leb ( B ( 0 , \delta ) )$
\end{itemize}
where $C (\epsilon) \ra 1$ when $\epsilon$ goes to zero.
\end{lem}
\begin{proof}
  By the hypothesis on the $graph ( h )$ and the definition of angle,
  we conclude that $ \| D_x h \| \leq \epsilon $ and the proof of the first
   item goes just by the Mean Value Theorem. The second item is also easy
   to prove just by the formula of volume for graph of a function
   (see Chapter 1 of \cite{Ca92} for the formulas).
 \end{proof}
In what follows we consider a $C^1$ function which is defined on a
ball of a linear subspace of $\mathbb{R}^n$ to another subspace.
We show the relation between the norm of the derivative of such
function and another one which locally has the same graph and is
defined on a slightly perturbed domain or codomain.
\begin{lem} ~\label{moft}
If $h$ is a $C^1$ function from  $B ( 0 , r ) \subset E $ to $ F $
such that
 $ \| D h ( x ) \| \leq a $ (small)  where $ F $ is a linear subspace
  with $ \measuredangle ( E^{\perp} , F ) \leq b $ (also small) then
   $ graph ( h ) $ will be graphic of a new function
   $ \tilde{h} : Dom ( \tilde{h}) \subset E \rightarrow  E^{\perp}$ and
   $ \| D \tilde{h} \| \leq K a  $ (where $ K $ is constant
   converging to 1 when $b$ goes to zero)
\end{lem}
\begin{lem} ~\label{desvio}
 If $h$ is a $C^1$ function from  $B ( 0 , r ) \subset E $ to $ E^{\perp}$ such that
  $ \| D h ( x ) \| \leq a $ (small) and $ F $ is a linear subspace with
   the same dimension of $E$ , with $ \measuredangle ( E , F ) \leq b $
    (also small) then $ graph ( h ) $ will be graphic of a new function $ \tilde{h} : Dom ( \tilde{h}) \subset F \rightarrow  F^{\perp}$ and $ \| D \tilde{h} \| \leq 2 ( a + b  ) $
\end{lem}
 The proof of the Lemma~\ref{moft} comes out just by the definition of angle and the
  derivative of a function. We prove Lemma \ref{desvio} as follows.
\begin{proof}
 First observe that by the definition of angle in Definition \ref{angel}, $ \| D h ( x ) \|$
 is equal to the angel between $E$ and $T_{( x , h ( x ) )} graph ( h )$. So,
 to prove the Lemma suppose that
 $ \measuredangle ( E , F ) = b$ and $ \measuredangle ( E , G ) = a $  with $ a , b $ small.
   Let $f$ be linear maps from $F$ to $F^{\perp}$ whose graph is $E$ 
   and $\tilde g : E \ra F^{\perp}$ and $g : F \ra F^{\perp}$ be the maps
    with $G$ as their graph. We are going to show that $\| g \| \leq 2 ( a + b )$.
    By the definition of
$\measuredangle ( F , E )$ and using Lemma \ref{moft} we have (see figure \ref{desv}):
\begin{figure}[t]
\psfrag{E}{$E$}
\psfrag{F}{$F$}
\psfrag{G}{$ G $}
\psfrag{E'}{$ E^{\perp} $}
\psfrag{F'}{$ F^{\perp} $}
\psfrag{x}{$ x $}
\epsfysize=2in
\centerline{\epsfbox{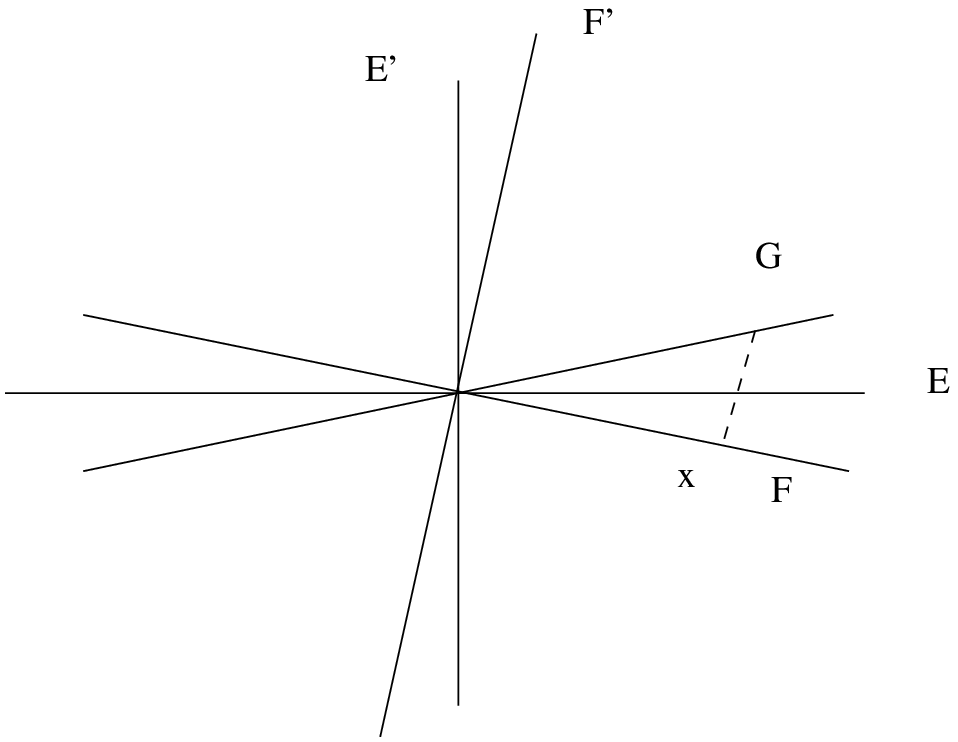}}
\label{desv}
\end{figure}

$$ \| g ( x ) \| \leq \| f ( x ) \| + \| \tilde g ( x + f ( x ) ) \| \leq \| f ( x ) \| + K a \| x + f ( x ) \|$$ where $K$ is near to one and is obtained by Lemma \ref{moft}. so we get :
$$ \frac{\| g ( x ) \|}{ \| x \|} \leq b + K a ( \sqrt{ 1 + b^2} ) \leq 2 ( a + b )$$ and the proof of Lemma \ref{desvio} is complete just by taking
$G = T_{ ( x , h ( x ) ) } graph ( h )$.
\end{proof}
The dependence of invariant subbundles $E^{cu}, E^{cs} $ to the base point is
 an important staff for the proof of the Theorem \ref{absolute}.
 The following control of the angles is a product of dominated decomposition and can be done with the same arguments as in \cite{Sh87}, pages 45-46. 
%
%
%
%
%
\begin{lem} ~\label{Mart}
There exist constants $ 0 < \alpha < 1 $ and $ 0 < \theta < 1$
with following property:

if $ d ( f^i ( x ) , f^i ( y ) $ is small for $ i = 0 , ... , n $ then for any two subspaces $ S_1 , S_2 $ respectively in $ C^{ cu } ( x ) , C^{ cu } ( y ) $ (small cones):
$$ \measuredangle ( D f_x ^n ( S_1 ) ,  D f_y ^n ( S_2 ) ) \preceq \theta ^n  +  dist (  f^n ( x )  , f^n ( y  ) )^{\alpha } . $$
\end{lem}
\begin{rem}
In the above Lemma, $\theta < 1 $ comes from dominated splitting and we can take
 $\theta^2 = \lambda$ where $\lambda$ is as in Lemma \ref{dominada}.
\end{rem}
%
%
%
%
%
%
\subsection {Covering $f^n ( S' )$ by graph of $C^1$ functions} \label{covering}
 We are going to show that for every point $ x $ in $ S' \cap D $,
$ f^n ( D ) $ locally can be seen as graph of a $ C^1 $ function
 from $ E^{cu} ( f^n ( x ) ) $ to $ E^{cs} ( f^n ( x ) ) $ with norm
  of derivative converging to zero uniformly as $ n $ goes to infinity.
   By this we intend to cover $f^n(S') \cap f^n (D)$ by flat disks.
    Let us call $ y_n := f^n ( x ) , y_n' := \pi_n ( y_n ). $

 We mention that for all $n$, $f^n ( D )$ is tangent to a thin cone which  varies continuously.  We show that there is a disk (inside $f^n(D)$) around $y_n$ which can be described as the graph of a $C^1$ function. The size of this disc decays when $n$ grows up, but it is definitely larger than the stable distance ($d_s$) between $y_n$ and $y_n'$.

%
%
%
%

\begin{lem} ~\label{kere}
There exists  $ \delta > 0 $ such that for $ \delta_1 < \delta $ and any $ x \in M $ If
$$ h : B_{\delta_1}^{cu} ( 0 ) \subset E^{cu} ( x ) \ra E^{cs} ( x ) ,\, h ( 0 ) = 0 ,\,  \| D h ( \xi ) \| \leq k ,\, \forall \xi \in B_{\delta_1}^{cu} $$
and
graph $( h ) \subset B_{\delta}^{cu} \times  B_{\delta}^{s}$ , then

$ W = f ( graph ( h ) ) \cap (  B_{ \gamma \delta_1}^{cu}  \times B_{\delta}^{s}  ) $  will be also graph of some $ \tilde h $ with the following properties:
\begin{enumerate}
\item Its domain contains $ B_{ \gamma \delta_1}^{cu}  $ and
                   $ \tilde h ( 0 ) = 0 ;$
\item  $\| D \tilde h ( \tilde{\xi}) \| \leq
                       k \theta , \, \forall \tilde{\xi} \in  B_{ \gamma \delta_1}^{cu}  \subset E^{cu} ( f ( x ) );$
\item $ {\bar \lambda} < \gamma $ where ${\bar \lambda}$ is defined as
(\ref{landabar}) in
 Section \ref{non-uniform}.
\end{enumerate}
  \end{lem}
\begin{proof}
 As $ f $ is $ C^2 $ , there exists $ \delta $ such that for all $ x \in M $, $ f$ can be written as :
 $ f ( \xi , \eta ) = ( A^{cu} ( \xi) + \phi^{cu} ( \xi , \eta ) , A^{cs} ( \eta ) + \phi^{cs} ( \xi , \eta ) $ , where $( \xi , \eta ) \in  B_{\delta}^{cu} \times  B_{\delta}^{cs} $ and $ \| D (  \phi^{cu} , \phi^{cs} ) \| \leq \epsilon $ . Just to reduce the notations suppose that $x$ is a fixed point. We define :
\begin{itemize}
 \item $ \alpha ( \xi ) = \tilde \xi: =   A^{cu} ( \xi) + \phi^{cu} ( \xi , h ( \xi ) = A^{cu}  ( \xi + (A^{cu})^{-1} ( \xi) \phi^{cu} ( \xi , h ( \xi ) ) $;
 \item $ \beta ( \xi ) =  A^{cs} ( h ( \xi ) ) + \phi^{cs} ( \xi , h ( \xi) ) ;\hspace{.5 in} \xi \in  B_{\delta}^{cu} ( 0 ); $
\end{itemize}
Now as $ \| ( A^{cu} )^{-1} \| \leq 1 + \delta_0 $ choosing $\epsilon$ small enough we deduce that:
              $$ \| ( A^{cu} )^{ - 1 } \| Lip ( \phi^{cu}
                          ( \xi , h ( \xi ) ) )  <  1 , $$
and this shows that $ \alpha = A^{cu} ( . ) ( Id +  (A^{cu})^{-1} (.) \phi^{cu} ( . , h ( . )) $ is invertible. So, it is enough to determine the domain of $ \alpha^ { - 1 }$ and defining $ \tilde h = \beta \circ \alpha^{ - 1 } $ for proving the first part of the lemma.\\
Observe that
                $$ \| \alpha ( \xi ) \| \geq \| A^{cu} ( \xi) \| -
                      \| \phi^{cu} ( \xi , h ( \xi ) \|
                 \geq ( \frac{1}{ 1 + \delta_0} - 2 \epsilon ) \| \xi \| >                                        \gamma \, \| \xi \| \,\,\, ,$$
where $\gamma$ is near to one as $\delta_0$ is small enough.
 Now by the aid of the proof of the inverse function theorem $ \alpha^{ - 1 } $ is defined on $ B_{ \gamma \delta_1}^{cu} $ and $$\tilde h = \beta \circ \alpha^{ - 1 } : B^{cu}_{\gamma \delta_1} \ra B^{cs} $$ is what we want.
  Observe that as ${\bar \lambda} < 1$, the third part of the lemma also turns out.

 Now we will verify the claim about derivative of $\tilde h $.
 By dominated splitting we have $ 0 < \theta < 1$ such that
      $ \| (A^{cu})^{-1} ( f ( x ) ) \| \, \| A^{cs} ( x ) \|
        \leq \theta^2$
($\theta^2$ is just the $\lambda$ in Lemma \ref{dominada}).
 By choosing $\epsilon $ small enough such that
  $ \| D \beta \| \leq \frac{\displaystyle k}{\displaystyle \sqrt\theta}$
  we get
            $$ \| D \tilde{h} ( \tilde{\xi} ) \| \leq \| D \beta (\xi)\| \,
            \| D \alpha^{-1} ( \tilde{\xi} ) \| \leq \frac{k}{\sqrt\theta} \|
               A^{cs} \| \, \| (A^{cu})^{-1} \| \, \| D ( I + T )^{-1}\|  $$
where $T = (A^{cu})^{-1} \phi^{cu} ( \xi , h ( \xi ) ) $ , on the other hand we have $$ \| D ( I + T )^{-1} \| = \| ( I + DT )^{-1}\| \leq \sum_{i=0}^{\infty} \| (DT)^i \| = \frac{1}{1 - \|DT\| } \leq \frac{1}{\sqrt\theta}$$
for $\epsilon$ is small enough. so $ \| D \tilde{h} (x)\| \leq k \theta $
\end{proof}
%
%
%
%
%
Let us see how to cover $f^n ( S ' ) \cap f^n ( D )$ by disks:

For $x \in S' \cap D$ there exists $ \delta > 0  $ (uniform in $D$) and $C^1$ functions $ h_x $
such that $ h_x : E^{cu}_{\delta} ( x ) \ra    E^{cs} ( x )$, and
the graph of $ h_x $ is a ball around $ x $. Now by Lemma \ref{kere} there exists $h_{ f^n
( x )} : B^{cu}_{ {\gamma}^n \delta} ( x ) \ra    E^{cs} ( f^n ( x
) ) $  such that $  h_{ f^n ( x )} ( B^{cu}_{ {\gamma}^n \delta} (
x ) )$ is a ball around $ y_n $ and also we have a good control on
the derivative of them : $ \| Dh_n ( x ) \| \leq k \, \theta^n$
where $h_n$ represents any $h_{f^n(x)}$. Applying Lemma
\ref{fasele} we get:
           $$ d ( z , y_n ) \leq d_1 ( z , y_n ) \leq
          k_n d ( z , y_n )  \quad \text{$\forall  z \in$ graph$( h_n )$ and  $k_n \ra 1$ } \quad ,
          $$
 and this gives that $ h_{ f^n ( x )} ( B^{cu}_{ \gamma\delta} ( x ) )$
is a ball of radius arbitrary near to $ 2a_n := \gamma^n \delta $  by taking $n$
large enough.
 we call this ball $ {\bar B_n }$ (around $ y_n $)  and $ B_n $ the ball with radius $ a_n $ around $ y_n $.

we mention that $\bar B_n $ is also graph of a function from
$E^{cu} $ to $ (E^{cu})^{\perp}$ over $ P ( \bar B_n  )$ where $P
$ is the orthogonal projection along $ (E^{cu})^{\perp}$.
\begin{rem} \label{dominio}
 By the estimate of the derivative of $h_n$,  $ P ( \bar B_n  )$ is contained in the ball of radius $ 2a_n ( 1 + C \theta^n )$ and contains  the ball of radius $ 2a_n ( 1 - C \theta^n )$ where $ C $ depends on the angle of $(E^{cu})^{\perp}$ and $E^{cs}$.
\end{rem}
In what follows we are working with $\bar B_n$ as the graph of the
mentioned new $ C^1$ function  which we call it also $h_n$ and it
is easy to see that $\| Dh_n \| \leq K \theta^n $ (Lemma
\ref{moft}).

Now we define a new transformation from $ {\bar B_n }$ to $ W^u ( q ) $ which is very near to holonomy $ \pi_n $.
Let's $ z \in  {\bar B_n }$ and define $ \cP ( z ) $ by translation along $ E^{cu}( y_n )^{\perp} $ which is orthogonal to the tangent space of all points of ${\bar B_n } $. One important property of $ \cP $ is that $ d ( z , \cP ( z ) ) $ is exponentially small. Indeed, we choose $a_n$ small enough for $ d ( z , \cP ( z ) ) $ being comparable to the $ d ( y_n , y_n' )$ = $ {\bar \lambda}^n $.
%
%
%
%
%
%
\subsection {Comparing measures of $B_i$ and $\pi_n ( B_i )$:} \label{comparing}
In the previous section we saw how to cover $f^n ( S ' ) \cap f^n
( D )$ by balls $\bar B_i$. In what follows we prove that the
volume of these disks does not increase ``a lot" by holonomy.
Indeed, we have to take $a_n$ in a good way to have this property.
The most important property for $a_n$ is:
\begin{equation} \label{hamishe}
 \frac{\displaystyle {\bar \lambda}^n}{\displaystyle a_n}  \ra 0
\end{equation}
 and the main proposition is the following.
\begin{pro} ~\label{compare}
There is a constant $I > 0$, independent of $n$, such that $ Leb ( \pi_n ( B_i ) ) \leq I Leb ( B_i )$
\end{pro}
To prove the above Proposition, we start with some lemmas.
%
%
%
%
\begin{lem} \label{hendese}
There is a choice of $a_n$ satisfying (\ref{hamishe}) such that for every $ z \in \bar{B_i}, d ( z , \cP ( z ) ) \preceq  {\Bar \lambda}^n $.
\end{lem}
\begin{proof}
 When $n$ is large enough we can consider $\cP ( \bar{B_n} ) $ also as a graph over $E^{cu} ( y_n ) $ to $ E^{cu} ( y_n )^{\perp}$, but we have to consider the angle between $E^{cu} ( y_n ) $ and $ E^{cu} ( y_n' )$ to calculate norm of derivative of the new function. To estimate norm of the derivative of the $C^1$ functions whose graphs are $\bar{B_n}$ and $\cP ( \bar{B_n} )$  we use lemmas \ref{desvio} and \ref{moft}.
 Using Mean Value Theorem and Remark \ref{dominio} we have (see figure \ref{hend}) :
\begin{figure}[t]
\psfrag{Ecu}{$E^{cu} ( y_n )$}
\psfrag{Ecuyn}{$ E^{cu} ( y_n )$}
\psfrag{Ecuyn,}{$ E^{cu} ( y_n' )$}
\psfrag{yn}{$ y_n$}
\psfrag{yn,}{$ y_n'$}
\psfrag{Ecuo}{$ E^{cu} \perp $}
\psfrag{Ecsyn}{$ E^{cs} ( y_n ) $}
\psfrag{Ecsyn,}{$ E^{cs} ( y_n' ) $}
\psfrag{pz}{$ \cP ( z )$}
\psfrag{z}{$ z $}
\epsfysize=4in
\centerline{\epsfbox{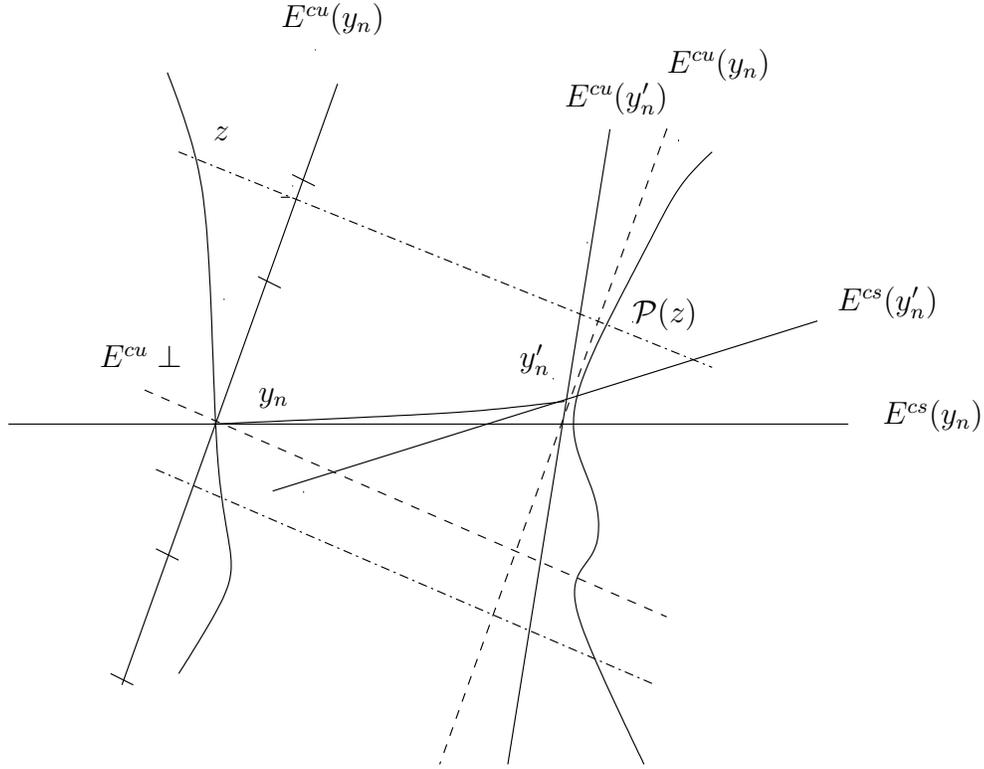}}
\caption{ Image of graphics}
\label{hend}
\end{figure}
      $$ d ( z , \cP ( z ) ) \leq  K ( 2 a_n + 2 C a_n \theta^n ) \theta^n +$$
      $$ +  ( 2 a_n + 2 C a_n \theta^n ) ( K \theta^n + \measuredangle ( E^{cu} ( y_n ) ,  E^{cu} ( y_n' ) ) )+ \bar{\lambda}^n   $$
Note that the term containing angles, in the above relations is
because of the deviation of $E^{cu} ( y_n ) $ from $ E^{cu} ( y_n'
)$ and applying lemma~\ref{desvio} :
                $$\measuredangle ( E^{cu} ( y_n ) ,  E^{cu} ( y_n' ) ) \preceq
                \theta^n + d ( y_n , y_n' )^{\alpha} , $$
and so:
       $$ d ( z , \cP ( z ) ) \preceq  a_n \theta^n +  a_n  d ( y_n , y_n' )^{\alpha}  + \bar{\lambda}^n . $$
So, to finish the proof of Lemma \ref{hendese} it is enough to choose $a_n$ satisfying
 the following two conditions : \\
\begin{itemize}
\item $ a_n \approx {\Bar \lambda} ^n \theta^{ - n }$
\item $ a_n \approx {\Bar \lambda}^{ n ( 1 - \alpha ) } $
\end{itemize}
Remember that by Lemma \ref{kere}, we need another restriction on
$ a_n $ to have graph of functions to use Mean Value Theorem.
\begin{itemize}
\item $ a_n \leq \gamma^{n} \delta $
\end{itemize}
 So choose $ a_n = min ( { \Bar \lambda} ^n \theta^{ - n } , { \Bar \lambda}^{ n ( 1 - \alpha ) } , {\gamma}^n \delta )$. As $\bar \lambda < \gamma $, already
$\frac{\displaystyle \bar{\lambda}^n}{\displaystyle a_n} \ra \, 0.$
\end{proof}
\begin{lem} \label{kocho} $\pi_n ( B_i )$ is contained in a ball around $ y_n'$ of radius near enough to $ \frac{\displaystyle 3}{\displaystyle 2} a_n $ as $n$ is large enough.
\end{lem}
\begin{proof}
For $z \in B_i $, $ \pi_n (z)$ lies in the $W^u ( q ) $ which is
contained in the graph of a function defined globally and the
graph is tangent to a thin cone field. So, by Lemma~\ref{fasele}
we deduce that
 for $z \in B_i $,
  $d_2 ( \pi_n ( z ) , y_n' ) \leq  \frac{\displaystyle 3}{\displaystyle 2}  d ( \pi_n ( z ) , y_n' ) $
so we get :
  $$ d_2 ( \pi_n ( z ) , y_n' ) \leq \frac{3}{2}  ( d ( \pi_n ( z ) , y_n' ) ) \leq \frac{3}{2} \, ( d ( \pi_n ( z ) , z ) + d ( z , y_n ) + d ( y_n , y_n' ) )$$
  $$ \leq \frac{3}{2} ( d_s ( \pi_n ( z ) , z ) + d_1( z , y_n ) + d_s( y_n , y_n' )  \leq \frac{3}{2} ( {\Bar \lambda}^n + k_n a_n + {\Bar \lambda}^n ) $$
$$ = \frac{3}{2} a_n ( k_n + \frac{\displaystyle 2{\Bar \lambda}^n }{\displaystyle a_n} )$$
So, choosing $a_n$ as in the Lemma~\ref{hendese} the proof is complete.
\end{proof}
\begin{lem} ~\label{gonde}
$ \cP ( { \bar B_i } ) $ contains $\pi_n ( B_i )$.
 \end{lem}
\begin{proof}
 For every $ z \in {\bar B_i} $ by triangular inequality for distances on the ambient manifold :
$$ d_2 ( \cP ( z ) , y_n' ) \geq d ( \cP ( z ) , y_n' ) \geq d ( z , y_n ) - d ( y_n' , y_n ) - d ( z , \cP ( z ) ) $$
$$ \geq \frac{1}{k_n} d_1( z , y_n ) - d_s( y_n' , y_n ) - d ( z , \cP ( z ) ) \geq  \frac{1}{k_n} d_1  ( z , y_n ) - 2 { \bar \lambda}^n .$$
 Here we have used Lemma \ref{hendese}, as
$ k_n \ra 1 $ and
 $ \frac{\displaystyle {\bar \lambda}^n}{\displaystyle a_n}  \ra 0 $
we conclude that $ \cP ( \bar B_i ) $ contains a ball around with
radius near to $ 2 a_n$ and by Lemma \ref{kocho} it contains $
\pi_n ( B_i )$.
\end{proof}
\begin{proof} (of Proposition \ref{compare}) Choose $a_n$ as in Lemma \ref{hendese}.
%
%
%
%
 As $Leb ( \bar B_i ) \leq I_1 Leb ( B_i ) $ , for a constant $ I_1 $ not depending to $n $
 and just depends to dimension of  $B_i$, we have :
$$ Leb ( \pi_n ( B_i ) ) \leq Leb ( \cP ( \bar{B_i} ) ) \approx Leb  (\bar{B_i} )
\leq I_1 Leb ( B_i ) $$ and the proposition is proved.
\end{proof}
Up to now we have covered $S_n:= f^n ( S' ) \cap f^n ( D ) $
 by a family of disks such that the volume of whose images under holonomy
  is comparable to their volume.
  By Besicovich covering theorem
  ~\cite{Ma95} we can cover $S_n$ with a countable locally finite
   subfamily $\{ B_i \}_i$ that is, there is a constant $C$ only
    depending to the dimension of $D$ such that, the intersection
     of any $C + 1$ disk of such subfamily is empty set.\\
%
%
%
%
\subsection{ Distortion estimates}\label{distortion}

 Now we state the distortion controls statements. By $ J f ( x , A ) $ we mean $ det ( D f _x | A )$
\begin{lem} [Bounded Distortion]
 There are $P_1 , M > 0$ such that for any $ z \in B_i $ the followings are satisfied:
\begin{itemize}
\item $ \frac{\displaystyle 1}{\displaystyle M} \leq \frac{\displaystyle J f^{ - n } ( y_n ,  T_{y_n } B_i  )}{ \displaystyle J f^{ - n } ( y_n' , T_{y_n' } \cP ( B_i )  )} \leq  M   $
\item
$ \frac{\displaystyle 1}{ \displaystyle P_1} \leq \frac{ \displaystyle J  f^{ - n } ( z , T_z ( B_i ) )}{ \displaystyle J f^{ - n } ( y_n , E^{cu}_{ y_n } )} \leq P_1 $
\end{itemize}
\end{lem}
\begin{proof} The problem is that in general we do not have H\"older control of the centre unstable fibers.
But in the case of the dominated decomposition or in other words when we have hyperbolicity property for the angles, one can show statements near to H\"older continuity .

As $ f $ is $ C^2$  function, we conclude that there exist
constants $ R_1 , R_2  > 0$ such that if $ z_1 , z_2 \in M  , d (
z_1 , z_2 ) \leq 1 $ and $ S_1 , S_2 $ are subspaces of $
\mathbb{R}^n $ with dimension  $u$ (dimension of $E^{cu}$) then:
\begin{equation} \label{hoho}
  |  \log J  f^ { -1 } ( z_1 , A_1 ) -  \log J  f^ { -1 } ( z_2 , A_2 ) | \leq R_1  d ( z_1 , z_2 )  + R_2 \measuredangle ( A_1 , A_2 ) .
\end{equation}

Now using the above inequality and Lemma \ref{Mart} we have :
 $$ | \log J  f^ { -n } ( y_n  , E^{ cu } ( y_n ) ) - \log J  f^ { -n } ( y_n'  , E^{ cu } ( y_n' ) ) | \leq $$
 $$ \leq R_1 ( \sum_{ i = 0 }^{ n - 1 } dist (  f^{ - i } ( y_n )  , f^{ - i } ( y_n'  ) ) ) + R_2 ( \sum_{i = 0 }^{n - 1 }  \measuredangle ( E^{cu} ( f^{- i } ( y_n ) , E^{cu} ( f^{- i } ( y_n' ) ) ) )$$

$$ \leq \frac{ C R_2}{ 1 - \theta } + ( K R_2 + R_1 ) \sum_{ i = 0 }^{ n - 1 } dist (  f^{ - i } ( y_n )  , f^{ - i } ( y_n'))^{\alpha }. $$

for some constants $C, K > 0$. So, using another time (\ref{hoho}) we conclude :
 $$ | \log J  f^ { -n } ( y_n  , T_{y_n } B_i  ) -  \log J  f^ { -n } ( y_n ' , T_{y_n' } \cP ( B_i ) ) | \leq      $$

$$  | \log J  f^ { -n } ( y_n  , T_{y_n } B_i  ) - \log J  f^ { -n } ( y_n  , ( E^{ cu } ( y_n  ) ) ) | +  $$

$$ | \log J  f^ { -n } ( y_n  ,  E^{ cu } ( y_n  )  ) - \log J  f^ { -n } ( y_n'  ,  E^{ cu } ( y_n' ) ) | + $$

$$| \log J  f^ { -n } ( y_n'  ,  E^{ cu } ( y_n'  )  ) - \log J  f^ { -n } ( y_n'  ,  T_{y_n' } \cP ( B_i ) ) |  \leq $$
\begin{equation} \label{haha}
 \frac{R_2}{ 1 - \theta } + ( K R_2 + R_1 ) \sum_{ i = 0 }^{ n - 1 } dist (  f^{ - i } ( y_n )  , f^{ - i } ( y_n' ) )^{\alpha } + 2 R_2 \sum_{ i = 0 }^{ n -1 } \theta^{ n - i }
\end{equation}
 As  $ y_n , y_n'$ are on the same strong stable manifold all of the terms appeared in (\ref{haha}) are summable and the proof of the first item of the lemma is complete. In fact, our argument show that we can substitute $y_n, y_n'$  respectively by any point $w_n \in B_i \cap f^n ( S' )$ and $\pi_n ( w_n )$.    The second item of the lemma comes out from the same arguments remembering that the size of $B_i \subset f^n ( D )$ is exponentially small.
\end{proof}
Now we apply distortion estimates of jacobians to get $$ Leb ( \pi ( D ) ) \leq \sum_i Leb ( f^{ - n } ( \pi_n ( B_i ) ) \leq   M P_1^2 \sum_i Leb ( f^{ - n } ( B_i )) \frac{ Leb ( \pi_n ( B_i ) ) }{ Leb ( B_i ) } $$
$$ \leq I M P \sum_i Leb ( f^{- n } ( B_i ) ) $$

But as $\{ B_i\}_i$ is a locally finite family covering $S_n$ and by $f^{-n}$ the areas of disks tangent to $C^{cu}$ decreases, taking $ n $ sufficiently large we see that
$\sum_i Leb ( f^{ - n } ( B_i ) ) \leq \, A Leb  ( D ) $. So taking $P_1^2 = P$ we conclude

$$ \frac{Leb ( \pi ( D ) )}{ Leb  ( D )} \leq IMPA ( universal ) $$
\newpage
\section {Appendix A: Robust indecomposability}
Topological transitivity of $C^1$ diffeomorphisms and ergodicity
(metric transitivity) of the Lebesgue measure for the $C^2$
conservative systems are two kinds of indecomposability. The
existence of SRB measures with full support and full Lebesgue
measure of basin (like in $C^2$-Anosov diffeomorphisms case) is
also a kind of indecomposability which in conservative diffeomorphisms case
 implies ergodicity. By results of \cite{BDP} we know that
$C^1$-robust transitivity implies dominated splitting. On the
other side, the results in robust ergodicity are for $C^2$
diffeomorphisms. For constructing SRB measures we need also more
regularity than $C^1$. So, we define $C^1$-robust
indecomposability as following:
\begin{defi} \label{indecomp}
Let $\diff^{1+} = \cup_{\alpha > 0} \diff^{1 + \alpha} (M)$. For
$f \in \diff^{1 +}$ we say $f$ is $C^1$-robustly indecomposable if
there is an open set $U \subset \diff^1(M)$ such that any $g \in U
\cap \diff^{1 +}$ has an SRB measure with $\mu$ $Leb ( B(\mu)) =
1$ and $Supp (\mu) = M$.
\end{defi}
\begin{pro}
Any $C^1$-robustly indecomposable diffeomorphism has dominated
splitting.
\end{pro}
\begin{proof}
Let $U$ be an open set as in the Definition\ref{indecomp}. We
claim that any $f \in U \cap \diff^{1+} (M)$ is transitive. To
show this, take two open sets $A, B$ in $M$. As $Supp (\mu) = M$ so,
$\mu (A), \mu (B) > 0$. Let $x \in B (\mu)$, by definition of the
basin, the orbit of $x$ goes through $A$ and $B$ infinitely many
times. This means that some iterate of $A$ intersects $B$.

 Now suppose $g_1 \in U$
does not admit dominated splitting, by the results in \cite{BDP}
one can perturb $g_1$ to get $g_2 \in U$ with a sink. Now by
density of $\diff^{1+} (M)$ in $\diff^1 (M)$ and persistence of
sinks in $C^1$ topology we get a diffeomorphism $g_3$ in $U \cap
\diff^{1+} (M)$ which has a sink and so can not be transitive
contradicting the above claim.
\end{proof}
However in the conservative case, the similar question is open.
\begin{que}
Does $C^2$ robust ergodicity or even $C^1$ robust ergodicity
defined as definition \ref{indecomp} imply dominated splitting.
\end{que}
Very roughly speaking by these results and questions we would like
to state: ``A robust indecomposability for dynamical systems requires some
weak form of hyperbolicity".

\section {Appendix B: Simultaneous hyperbolic times}
In ~\cite[Theorem 6.3]{ABV00} ergodic $cu-$Gibbs measures for
diffeomorphisms with dominated splitting and the non-uniformly
hyperbolic property like in the Preliminary section, are
constructed. This measures are absolutely continuous along a
family of disks which are tangent to center-unstable conefield.
\begin{pro}
For $f \in \cV$, the $cu-$Gibbs measures as above are SRB i.e
their basin has positive Lebesgue measure.
\end{pro}

 To prove that these measures are SRB, one need to show that for points
  in the support of these measures, all the Lyapunov exponents
  (in the $E^{cs}$ direction) are negative. To provide negative
  Lyapunov exponents, in ~\cite{ABV00}, the authors add the
  condition of ``simultaneous hyperbolic times''. We show that
   for $f \in \cV$ it is not necessary to verify this condition
   and see that the $cu-$Gibbs measure constructed there, are indeed SRB measure.

 For any $y \in Supp (\mu) $ where $\mu$ is one of such $cu$-Gibbs measures, there exists $x$ such that $y \in D^{\infty} ( x )$ where $D^{\infty} ( x )$ is tangent $E^{cu}$ at any point of it
  and moreover it is the local strong unstable manifold of $x$ (see ~\cite[Lemma 3.7]{ABV00}).
\begin{lem}
If $f \in \cV$ then for Lebesgue almost all point of  $D^{\infty}
( x ) $ the Lyapunov exponents in the $E^{cs}$ direction are
negative.
\end{lem}
\begin{proof}
By the above observations about $D^{\infty} ( x )$ we may consider
the lift of $D^{\infty} ( x ) $ to $\mathbb{R}^n$ included in the
graph of a global $C^1$ function $\gamma: \mathbb{R}^u \ra
\mathbb{R}^s$ with $T_{( z , \gamma ( z) )} graph (\gamma) \in
C^{cu} ( z , \gamma ( z) )$. So, by the definition of dynamically
flat submanifols in Section \ref{non-uniform}, $D^{\infty} ( x )
$ is contained in a dynamically flat submanifold and by
Proposition ~\ref{almost} for almost all points in $D^{\infty} ( x
) $ all the Lyapunov exponents in the $E^{cs}$ direction are
negative.
\end{proof}
For proving that the $cu-$Gibbs measures are really SRB, or the
basin of them has positive volume, we repeat the same argument of
~\cite[Proposition 6.4]{ABV00}:
\begin{proof}
Let $\mu$ be such a Gibbs ergodic measure. There exists some disk $D^{\infty}$ such that almost every point in $D^{\infty}$ is in the basin of $\mu$. By absolute continuity of stable lamination of the points in $D^{\infty} \cap B ( \mu )$ and the fact that these stable manifolds are contained in $B ( \mu )$, we conclude that the basin of $\mu$ must have positive Lebesgue measure.
\end{proof}
\bibliographystyle{alpha}
\bibliography{bib}
\end{document}